\documentclass[12pt]{article}

\usepackage{amsfonts}
\usepackage{amssymb}
\usepackage{amsmath}
\usepackage{amsthm}
\usepackage{epsfig}
\usepackage[mathscr]{eucal}
\def\RR{{\mathbb R}}

\def\CC{{\mathbb C}}

\def\W{{\mathbb W}}
\def\E{{\mathbb E}}

\newtheorem{theorem}{Theorem}
\newtheorem{pro}[theorem]{Proposition}

\newtheorem{lemma}{Lemma}
\newtheorem{corollary}{Corollary}
\newtheorem{remark}{Remark}

\def\NN{{\mathbb N}}
\def\ZZ{{\mathbb Z}}
\def\con{{\cal C}}
\def\He{{He}}
\def\nHe{{Ne}}
\def\dHe{{\tilde{He}}}
\def\dnHe{{\tilde{Ne}}}

\def\eop{$\quad \Box$}
\def\meop{\quad \Box}

\def\proof{{\bf Proof:} \,}

\begin{document}

\title{Interpolation of exponential-type functions
on a uniform grid by shifts of a basis function\\
\bigskip
{\small Dedicated to the memory of Ward Cheney}}

\author{A. K. Kushpel, J. Levesley, X. Sun\thanks{This work was supported by EPSRC Grant EP/H020071/1, the University of Leicester via study leave for the second author and Missouri State University through a generous travel grant for the second author.}}

\bigskip

\date{}
\maketitle
\begin{abstract}
In this paper,  we study the problem of interpolating a continuous function at $(n+1)$ equally-spaced points in the interval $[0,1]$, using shifts of a kernel on the $(1/n)$-spaced infinite grid. The archetypal example here is approximation using shifts of a Gaussian kernel. We present new results concerning interpolation of functions of exponential type, in particular, polynomials on the integer grid as a step en route to solve the general interpolation problem. For the Gaussian kernel we introduce a new class of polynomials, closely related to the probabilistic Hermite polynomials and show that evaluations of the  polynomials at the integer points provide the coefficients of the interpolants. Finally we give a closed formula for the Gaussian interpolant of a continuous function on a uniform grid in the unit interval (assuming knowledge of the discrete moments of the Gaussian).
\end{abstract}

\section{Introduction}

  In the mathematical literature pertaining to radial basis functions, there have been mainly two approaches to constructing interpolants with Gaussian kernels. The first involves interpolating a function $f \in \con(\RR)$ on the $h$-spaced grid $h \ZZ$ by an interpolant of the form
$$
\sum_{z \in \ZZ} \alpha_z \psi(z-x/h),
$$
where
$$
\psi(x)={\exp(-\|x\|^2/2) \over \sqrt{2 \pi}}.
$$
Analysis of this so-called cardinal approximation has been done in a series of papers of Baxter, Riemenschnieder and Sivakumar \cite{baxter,kush,sivakumar,riemenschneider}.

The second concentrates on interpolating a continuous function on a finite subset $Y$ of a compact interval (e.g $[0,1]$). Under this circumstance, the interpolant one seeks is of the form
$$
\sum_{y \in Y} \alpha_y \psi(y-x).
$$
There are multidimensional set-ups for both approaches. Modern mathematical literature abounds in developing error estimates for approximation schemes in this context. We refer readers to \cite{hamm,hangelbroeck,larrson,madych} and the references therein.

Approximation methods involving sparse-grid algorithms have been recently proven effective and efficient; see \cite{griebel}. Some sophisticated multi-level sparse grid kernel interpolation schemes have been constructed by authors of \cite{dong,mlski}.
We are currently motivated to develop sparse-grid algorithms for high-dimensional approximation
with the Gaussian kernel and derive error estimates for $C^k-$ functions with polynomial growth. However, there are several obstacles en route to achieving these goals. The main purpose of the current paper is to clear a few obstacles out of the way.  First and foremost, we face the problem of interpolating a function at the $(n+1)$ equally-spaced points $ih, i=0,1,\cdots,n$ with $h=1/n$, where $n \in \NN$. The approach we take here differs from those discussed in the above references. We first interpolate the given $(n+1)$ data by a degree $n$ polynomial, and then interpolate the polynomial by a radial basis function interpolant on $h \ZZ$.

Functions of exponential type\footnote{Some variants of functions of exponential type are also referred to as ``band limited functions" in the literature.}  are often utilized as a half-way house in deriving error estimates for Sobolev space functions; see \cite{NWW, NSWW}. As such, it is worthwhile to study the effect of the interpolation scheme when the target functions are of exponential types, and in particular, polynomials, which we anticipate to play a significant role in our future effort to obtain more nuanced error estimates for $C^k-$functions with polynomial growth.
Interestingly enough, we observe in the analytic number theory literature that interpolation goes the opposite way in the sense that functions of exponential type are employed to approximate the Gaussian and other useful radial basis functions; see \cite{carneiro1, carneiro2}. We hope that interactions of the two seemingly
inverse research tracks will create synergistic results.

The layout of the paper is as follows. In Section~\ref{general}, we will consider a general kernel $\psi$ and study the operator ${\cal T}$ induced by the Toeplitz matrix $\psi(j-k)$, $j, k \in \ZZ$. The action of ${\cal T}$ on an $ f \in \con (\RR)$ takes the form:
 \[
 \sum_{j \in \ZZ}  \psi(j-k) f(j), \quad k \in \ZZ.
 \]
 We show that the operator is one-to-one on the space of polynomials, which are the subject of interest in the main body of the current paper. In the Appendix we show furthermore that the operator is one-to-one on a larger class of functions of exponential type.
Moreover, we demonstrate that the coefficients of an interpolant can be written in closed form.  In Section~\ref{gaussian}, we investigate the special case in which the Gaussian kernel is employed and the target functions are polynomials. We introduce new classes of polynomials resembling the classical probabilistic Hermite polynomials, and derive closed formulas for the coefficients of a Gaussian interpolant in terms of these  polynomials.
  In Section~\ref{newton}, we show how to use interpolate general functions on a uniform grid in the unit interval via interpolation by polynomials. This is not a convergent approximation scheme. However, it can be used as part of a residual correction scheme as is described in \cite{mlski}.

\section{Interpolation with general kernels} \label{general}

Let us assume we have an infinitely differentiable, positive function $\psi$ such that all the ``discrete moments"  $M_k$:
\begin{eqnarray}
M_k & = & \sum_{j \in \ZZ} j^k \psi(j), \quad k=0,1,\cdots,N, \label{momdef}
\end{eqnarray}
are finite.
Let $p_k (x)= x^k$, for $0 \le k \le N$. We seek coefficients $a_{k,j}, \; j \in \ZZ$, such that
\begin{eqnarray}
I[p_k](x) & := & \sum_{j \in \ZZ} a_{k,j} \psi(x-j), \label{idef}
\end{eqnarray}
interpolates $p_k$ at all the integers. In other words
\[
I[p_k](\ell) = \ell^k, \quad 0 \le k \le N, \quad \ell \in \ZZ.
\]
 We will show that if coefficients $a_{k,j}$ (as polynomials in $j$ of degree $k$) exist, then they are unique. We then construct such   coefficients (of polynomial form). Thus $I[p_k]$ is the unique interpolant with coefficients of polynomial form. Furthermore, these coefficients are constructible in a recursive fashion.

In the Appendix at the end of the paper we extend results of existence and uniqueness to include functions of exponential type.

\begin{lemma} \label{unique}
Let $r_k$ be a polynomial of degree $k$ for some $0 \le k \le N$. Suppose that
$$
f(\ell) := \sum_{j \in \ZZ} r_k(j) \psi(j - \ell)=0, \quad \ell \in \ZZ.
$$
Then $r_k \equiv 0$.
\end{lemma}
\proof We proceed by induction. If $k=0$ then $r_k=c$ for some constant $c$. Then,
$$
c \sum_{j \in \ZZ}  \psi(j - \ell) = 0,
$$
and since $\psi$ is positive, $c=0$.

Assume that the result holds true for all polynomials of degree $<k$. Let $\Delta f (x) = f(x+1)-f(x)$, $x \in \RR$, be the forward difference operator. Then
\begin{eqnarray*}
\Delta f (\ell) & = & f(\ell+1)-f(\ell) \\
& = & \sum_{j \in \ZZ} r_k(j) \psi(j - \ell-1)-\sum_{j \in \ZZ} r_k(j) \psi(j - \ell) \\
& = & \sum_{j \in \ZZ} r_k(j+1) \psi(j - \ell)-\sum_{j \in \ZZ} r_k(j) \psi(j - \ell) \\
& = & \sum_{j \in \ZZ} (r_k(j+1)-r_k(j)) \psi(j - \ell) \\
& = & \sum_{j \in \ZZ} \Delta r_k(j)\psi(j - \ell) = 0.
\end{eqnarray*}
Note that $\Delta r_k$ is a polynomial of degree $k-1$, which satisfies $\Delta r_k(j)=0$, $j \in \ZZ$. By the induction hypothesis, we have
that $\Delta r_k = 0$. It follows that $r_k$ is a constant. We use the induction hypothesis once more to conclude that $r_k$ is identically zero. \eop

The result of Lemma \ref{unique} can be equivalently stated as follows. There is no nontrivial polynomial $p$,
$\deg(p) \le N$, such that
$$
 \sum_{j \in \ZZ} p(j)  \psi(j - \ell) = 0, \quad \ell \in \ZZ.
$$
In the sequel, we will use the above fact without further declaration.
\begin{lemma} \label{discpolyconv}
For $k=0,1,\cdots,N$,
\begin{eqnarray*}
\sum_{j \in \ZZ} j^{k} \psi(j-\ell) & = & \sum_{i=0}^k  {k \choose i} M_{k-i} \ell^{i}, \quad \ell \in \ZZ.
\end{eqnarray*}
\end{lemma}

\proof Substituting the expression for $M_{k-m}$ from (\ref{momdef}) we have
\begin{eqnarray*}
\sum_{i=0}^k  {k \choose i} M_{k-i} \ell^{i} & = & \sum_{i=0}^k  {k \choose i}  \left \{ \sum_{j \in \ZZ} j^{k-i} \psi(j) \right \} \ell^{i} \\
& = & \sum_{j \in \ZZ} \psi(j) \left \{  \sum_{i=0}^k  {k \choose i} j^{k-i} \ell^i \right \} \\
& = & \sum_{j \in \ZZ} \psi(j) (j+\ell)^k \\
& = & \sum_{j \in \ZZ} j^k \psi(j-\ell),
\end{eqnarray*}
by simply renumbering the sum. \eop

\smallskip

Let $B_{k,i}= M_{k-i} {k-1 \choose i-1}$, $k=1\cdots,N+1$, $i=1,\cdots,k$. Since
\begin{equation}
\sum_{j \in \ZZ} j^{k-1} \psi(j-\ell) =  \sum_{i=1}^{k} B_{i,k} \ell^{k-1}, \quad \ell=1,\cdots,N+1, \label{Amat}
\end{equation}
we can use backward elimination to see that
\begin{eqnarray}
\ell^{k-1} & = & \sum_{j \in \ZZ} \left \{ \sum_{i=1}^k A_{k,i} j^{i-1} \right \} \psi(j-\ell), \quad k=1,\cdots,N+1, \label{polyco}
\end{eqnarray}
for some numbers $A_{k,i}$, $k=1,\cdots,N+1, \; i=1, \cdots,k$. Define the polynomial $a_k(j)=\sum_{i=1}^k A_{k+1,i+1} j^{i-1}$. Then the coefficients $a_{k,j}$ defined in (\ref{idef}) satisfy $a_{k,j}=a_k(j)$.

 An immediate consequence of Lemma~\ref{unique} and (\ref{polyco}) is
 \begin{theorem}
 The coefficients $a_{k,j}$ of the interpolant $I[p_k]$ given in (\ref{idef}) are unique. Moreover, they are polynomials  (in $j$) of degree $k-1$.
\end{theorem}

\smallskip

Let
\begin{eqnarray}
q_k(x) & = & \sum_{m=0}^k M_{k-m} {k \choose m} x^{m}. \label{qdef}
\end{eqnarray}
In the next theorem we prove a generating function style relationship for the polynomials $q_k$ and the interpolating functions $I[p_k]$. As a corollary of this theorem we give a recursive formula for the errors between $I[p_k]$ and $p_k$.
\begin{theorem} \label{th:interp}
For $k =0,1,\cdots,N$,
\begin{eqnarray}
\sum_{j \in \ZZ} (j-x)^k \psi(j-x) & = & \sum_{i=0}^k {k \choose i} q_{i}(-x) I[p_{k-i}](x). \label{interp}
\end{eqnarray}
\end{theorem}

\proof If we substitute (\ref{polyco}) into Lemma~\ref{discpolyconv} we have
\begin{eqnarray*}
\sum_{j \in \ZZ} j^{k} \psi(j-\ell) & = & \sum_{i=0}^k  {k \choose i} M_{k-i} \sum_{j \in \ZZ} a_i(j) \psi(j-\ell), \quad \ell \in \ZZ,
\end{eqnarray*}
and by the linear independence of $\psi(\cdot - \ell)$, $\ell \in \ZZ$,
\begin{eqnarray}
j^{k} & = & \sum_{i=0}^k  {k \choose i} M_{k-i} a_i(j), \quad j \in \ZZ. \label{inverse}
\end{eqnarray}
We now expand the left hand side of \eqref{interp} to get
\begin{eqnarray*}
\sum_{j \in \ZZ} (j-x)^k \psi(j-x) & = & \sum_{j \in \ZZ} \sum_{i=0}^k {k \choose i} (-x)^i j^{k-i} \psi(j-x) \\
& = & \sum_{j \in \ZZ} \sum_{i=0}^k {k \choose i} (-x)^i \left ( \sum_{m=0}^{k-i}  {k-i \choose m} M_{k-i-m} a_m(j) \right ) \psi(j-x) \\
& = &  \sum_{i=0}^k {k \choose i} (-x)^i \sum_{m=0}^{k-i}  {k-i \choose m} M_{k-i-m} \left ( \sum_{j \in \ZZ} a_m(j) \psi(j-x) \right ) \\
& = & \sum_{i=0}^k {k \choose i} (-x)^i \sum_{m=0}^{k-i}  {k-i \choose m} M_{k-i-m} I[p_m](x),
\end{eqnarray*}
In deriving the above equations, we have used (\ref{idef}) and (\ref{inverse}). Making use of the formula
$$
{k \choose i} {k-i \choose m} = {k \choose m} {k-m \choose i},
$$
we reorder the final sum above as follows.
\begin{eqnarray*}
\sum_{j \in \ZZ} (j-x)^k \psi(j-x) & = & \sum_{m=0}^k I[p_m](x)\sum_{i=0}^{k-m} {k \choose i} {k-i \choose m} M_{k-m-i} (-x)^i \\
& = & \sum_{m=0}^k {k \choose m} I[p_m](x) \sum_{i=0}^{k-m} {k-m \choose i} M_{k-m-i} (-x)^i \\
& = &  \sum_{m=0}^k {k \choose m} I[p_m](x)q_{k-m}(-x).
\end{eqnarray*}
In deriving the last equation, we have used (\ref{qdef}). \eop

We define the error of interpolation $E_k(x):=I[p_{k}](x)-p_k(x)$, and
$$
\chi_k(x): = \sum_{j \in \ZZ} (j-x)^k \psi(j-x) - M_k, \quad k=0,1,\cdots,N.
$$
 Next, we obtain a recursive formula which, upon an appropriate inversion,  allows us to write the errors in terms of the functions $\chi_k$ and $q_k$.
\begin{corollary}
For $k = 0,1,\cdots,N$,
\begin{eqnarray*}
\chi_k(x) & = & \sum_{i=0}^k {k \choose i} q_{i}(-x) E_{k-i}(x).
\end{eqnarray*}
\end{corollary}
\proof Since $I[p_{k-i}](\ell)=p_{k-i}(\ell)$, $i=0,1,\cdots,k$, $\ell \in \ZZ$, by Theorem~\ref{th:interp}, we have for $k=0,1,\cdots,N$,
$$
\sum_{i=0}^k {k \choose i} q_{i}(\ell) \ell^k = \sum_{j \in \ZZ} (j-\ell)^k \psi(j-\ell) = M_k.
$$
Since a polynomial is uniquely determined by its values on the integers we have
\begin{eqnarray*}
\sum_{i=0}^k {k \choose i} q_{i}(-x) x^k & = & M_k.
\end{eqnarray*}
Subtracting this equation from (\ref{interp}) we see that
\begin{eqnarray*}
\chi_k(x) & = & \sum_{i=0}^k {k \choose i} q_{i}(-x) I [p_{k-i}](x) - M_k \\
& = &  \sum_{i=0}^k {k \choose i} q_{i}(-x) I [p_{k-i}](x) - \sum_{i=0}^k {k \choose i} q_{i}(-x) p_{k-i}(x) \\
& = & \sum_{i=0}^k {k \choose i} q_{i}(-x) (I [p_{k-i}](x) - p_{k-i}(x)) \\
& = & \sum_{i=0}^k {k \choose i} q_{i}(-x) E_{k-i}(x). \quad \meop
\end{eqnarray*}

\begin{remark}
Theorem~\ref{th:interp} gives a recursive formula for computing the interpolant of any polynomial, as long as one knows the interpolants for polynomials of lower degree. Likewise, the result of the corollary expresses the error between a polynomial and its interpolant in the same fashion. If the Gaussian kernel is employed to do interpolation, then we will have more interesting information to offer. We will study this topic in the next section.
\end{remark}


\section{The Gaussian Kernel} \label{gaussian}

In this section we study exclusively the case of $\psi(x)={\exp(-x^2/2) \over \sqrt{2 \pi}}$.
Pertinent to the contents of this paper will be the probabilistic Hermite polynomials $\He_k$, $k=0,1,\cdots$. These may be defined in a number of ways, but for us perhaps the most appropriate one is via Rodrigues formula:
\begin{eqnarray*}
\He_k (x) & := & \exp(x^2/2) {d^k \over dx^k} \exp(-x^2/2), \quad k=0,1,\cdots.
\end{eqnarray*}
We have the following explicit representation of these polynomials (see e.g. \cite{askey}):
\begin{eqnarray*}
\He_k(x) & = & \sum_{i=0}^{\lfloor k/2 \rfloor} (-1)^{k-i} {k! \over i! (k-2i)!} {x^{k-2i} \over 2^i}, \quad k=0,1,\cdots,
\end{eqnarray*}
where ${\lfloor x \rfloor}$ is the greatest integer less than or equal to $x$. The polynomial $\He_k$ has a close cousin that is often referred to as the probabilistic polynomial of negative variance:
\begin{eqnarray*}
\nHe_k(x) & = & \sum_{i=0}^{\lfloor k/2 \rfloor} {k! \over i! (k-2i)!} {x^{k-2i} \over 2^i}, \quad k=0,1,\cdots,
\end{eqnarray*}
which has the same coefficients in absolute value, but the coefficients are all positive.

The probabilistic polynomials of negative variance arise very naturally in this study as they are the result of the continuous convolution of the Gaussian with the polynomials of appropriate degree:
\begin{lemma} \label{lem:negherm}
Let $\psi(x)=\exp(-x^2/2)$. Then
\begin{eqnarray*}
\nHe_k(x) & = & \int_{-\infty}^\infty y^k \psi(y-x) dy \\
& = & \sum_{i=0}^{\lfloor k/2 \rfloor} {k \choose 2i} C_{2i} x^{k-2i},
\end{eqnarray*}
where
$$
C_k = \int_{-\infty}^\infty y^k \psi(y) dy = \left \{ \begin{array}{ll}
(k-1)!!, & k \; {\rm even}, \\
0, & k \; {\rm odd}.
\end{array}
\right.
$$
\end{lemma}
\proof It is well-known (see e.g. \cite{papoulis}) that
$$
\int_{-\infty}^\infty y^k \psi(y) dy = \left \{ \begin{array}{ll}
(k-1)!!, & k \; {\rm even}, \\
0, & k \; {\rm odd}.
\end{array}
\right.
$$
To prove the first equation, we make a simple change of variable $w=y-x$:
\begin{eqnarray*}
\int_{-\infty}^\infty y^k \psi(x-y) dy & = & \int_{-\infty}^\infty (w+x)^k \psi(w) dw \\
& = &  \sum_{i=0}^k {k \choose i} x^i \int_{-\infty}^\infty w^{k-i} \psi(w) dw \\
& = &  \sum_{i=0}^k  {k \choose i} C_{k-i} x^i \\
& = & \sum_{i=0}^k  {k \choose i} C_{i} x^{k-i} \\
& = & \sum_{i=0}^{\lfloor k/2 \rfloor} {k \choose i} C_{2i} x^{k-2i}.
\end{eqnarray*}
We remind readers that all the odd degree terms have zero coefficients. If we substitute the value for $C_{2i}$ we see that $\nHe_k$ is the probabilistic Hermite polynomial of negative variance:
$$
\nHe_k(x) = \sum_{i=0}^{\lfloor k/2 \rfloor} {k! \over i! (k-2i)!} {x^{k-2i} \over 2^i}. \meop
$$

A fascinating relationship between $\He_k$ and $\nHe_k$ is the so-called umbral composition (see \cite{askey}):
$$
\sum_{i=0}^{\lfloor k/2 \rfloor} {k! \over i! (k-2i)!} {\He_{k-2i}(x) \over 2^i} = \sum_{i=0}^{\lfloor k/2 \rfloor} (-1)^i {k! \over i! (k-2i)!} {\nHe_{k-2i}(x) \over 2^i} = x^k.
$$
Using Lemma~\ref{lem:negherm} and the second equation above we have
\begin{eqnarray}
x^k & = & \sum_{i=0}^{\lfloor k/2 \rfloor} (-1)^i {k! \over 2^i i! (k-2i)!} \int_{-\infty}^\infty y^{k-2i} \psi(x-y) dy \nonumber \\
& = & \int_{-\infty}^\infty \He_k(y) \psi(y-x) dy, \label{contint}
\end{eqnarray}
so that we can recover the monomials by integrating against the probabilistic Hermite polynomials. Of course, this gives us an idea of what will happen in the discrete case.

To do this,  we need an analogue of the probabilistic Hermite polynomial for the discrete case. We define
\begin{eqnarray}
\dHe_k(x) & = & \sum_{i=0}^{\lfloor k/2 \rfloor} (-1)^{i} {k \choose 2i} M_{2i} x^{k-2i}, \quad k=0,1,\cdots , \label{dhedef}
\end{eqnarray}
where $M_{2k}$ is the discrete moment as defined in (\ref{momdef}). We also let
$$
\dnHe_k = q_k, \quad k=0,1,\cdots,
$$
where $q_k$ is as defined in (\ref{qdef}). In other words
\begin{eqnarray}
\dnHe_k (x) & = & \sum_{i=0}^{\lfloor k/2 \rfloor} {k \choose 2i} M_{2i} x^{k-2i}. \label{dnhe}
\end{eqnarray}

 Equation (\ref{contint}) suggests that a closed formula for   the interpolant $I[p_k](x)$ resembles
$$
\sum_{j \in \ZZ} \dHe_k(j) \psi(j-x), \quad x \in \RR.
$$
In the next result, we  will show that this is indeed the case (up to a constant very close to 1) for $k=0,1,2,3$. For $k \ge 5$ we need to make some modifications to $\tilde H_k$ for a closed form.

\begin{lemma} \label{lem:corr}
For $k=0,1,\cdots$, and $\ell \in \ZZ$,
\begin{eqnarray*}
\sum_{j \in \ZZ} \dHe_k(z) \psi(j-\ell) & = & \sum_{i=0}^{\lfloor k/4 \rfloor} \left \{ \sum_{m=0}^{2i} (-1)^m {k \choose 2m}{k-2m \choose 4i-2m} M_{2m} M_{4i-2m} \right \} j^{k-4i}.
\end{eqnarray*}
\end{lemma}
\proof By Theorem~\ref{discpolyconv}, Equations (\ref{dhedef}) and (\ref{dnhe}), we have, for $\ell \in \ZZ$,
\begin{eqnarray*}
\sum_{j \in \ZZ} \dHe_k(z) \psi(j-\ell) & = & \sum_{j \in \ZZ} \left \{ \sum_{i=0}^{\lfloor k/2 \rfloor} (-1)^{i} {k \choose 2i} M_{2i} j^{k-2i} \right \} \psi(j-\ell) \\
& = & \sum_{i=0}^{\lfloor k/2 \rfloor} (-1)^{i} {k \choose 2i} M_{2i} \left \{ \sum_{j \in \ZZ} j^{k-2i} \psi(j-\ell) \right \} \\
& = & \sum_{i=0}^{\lfloor k/2 \rfloor} (-1)^{i} {k \choose 2i} M_{2i} \left \{ \sum_{m=0}^{\lfloor (k-2i)/2 \rfloor} (-1)^{m} {k-2i \choose 2m} M_{2m} \ell^{k-2i-2m} \right \}.
\end{eqnarray*}
 Rearranging we obtain
\begin{eqnarray*}
\sum_{j \in \ZZ} \dHe_k(z) \psi(j-\ell) & = & \sum_{i=0}^{\lfloor k/2 \rfloor} \left \{ \sum_{m=0}^{i} (-1)^m {k \choose 2m}{k-2m \choose 2i-2m} M_{2m} M_{2i-2m} \right \} \ell^{k-2i} \\
& = & \sum_{i=0}^{\lfloor k/4 \rfloor} \left \{ \sum_{m=0}^{2i} (-1)^m {k \choose 2m}{k-2m \choose 4i-2m} M_{2m} M_{4i-2m} \right \} \ell^{k-4i},
\end{eqnarray*}
where we have used the fact that, if $i$ is odd, then
$$
\sum_{j=0}^{i} (-1)^j {k \choose 2j}{k-2j \choose 2i-2j} M_{2j} M_{2i-2j} = 0,
$$
which is true because
\begin{eqnarray*}
{k \choose 2j}{k-2j \choose 2i-2j} & = & {k! \over (k-2j)! (2j)!} {(k-2j)! \over (k-2i)! (2i-2j)! } \\
& = & {k! \over (k-2i+2j)! (2i-2j)!} {(k-2i+2j)! \over (k-2i)! (2j)! } \\
& = & {k \choose 2(i-j)}{k-2(i-j) \choose 2j}. \quad \Box
\end{eqnarray*}

As we can see, for $k=0,1,2,3$, the above lemma gives an exact formula. Interestingly, suppose that we replace $M_k$ by $C_k$, $k=0,1,\cdots$, then the correction terms above are all zero. This is why we get the umbral composition formula for the probabilistic Hermite polynomials. For higher degrees we need to modify the polynomial in the summation for interpolation. To this end we introduce the polynomials $Q_k(x)$, which we define by
\begin{eqnarray}
Q_k & = & {1 \over M_0^2} \nHe_k, \quad k=0,1,2,3, \label{q1}
\end{eqnarray}
and for $k=4,5,\cdots$,
\begin{eqnarray}
Q_k  =  {1 \over M_0^2} \nHe_k-\sum_{i=1}^{\lfloor k/4 \rfloor} \left \{ \sum_{j=0}^{2i} (-1)^j {k \choose 2j}{k-2j \choose 4i-2j} M_{2j} M_{4i-2j} \right \} Q_{k-4i}. \label{q2}
\end{eqnarray}
Using Lemma~\ref{lem:corr},  we immediately get the main result of this section
\begin{theorem}
For $k=0,1,\cdots$, and $\psi$ the Gaussian, we have
$$
I[p_k](x) = \sum_{j \in \ZZ} Q_k(j) \psi(j-x).
$$
\end{theorem}

\section{Gaussian Interpolant on a $h$-spaced points in an interval} \label{newton}

In this section, we give a recipe for computation of a new Gaussian kernel interpolant to a function defined at equally spaced points $X=\{0, \frac1n,\frac2n, \cdots, 1\}$. The construction of this interpolant utilizes the full $(1/n)-$ spaced infinite grid. As such, it is different from most of Gaussian kernel interpolants constructed with conventional procedures. However, for all the practical computational purposes, only a small number of centres outside of the interval of interpolation are required. The rapid decay of the Gaussian kernel offsets the error incurred by dropping terms (shifts) of the interpolant far from the interpolation interval.

In this case we seek an interpolant of the form
$$
I_n [f] (x) = \sum_{j \in \ZZ} a_j^n(f) \psi(j-nx), \quad x \in [0,1],
$$
for $f \in \con([0,1])$. To do this we follow the following recipe:
\begin{description}
\item{1.} Interpolate $f$ on $X$ with a degree $n$ polynomial
$$
P_n [f](x) = \sum_{i=0}^n c_i(f) t_i(x),
$$
where $t_i$, $i=0,\cdots,n$ is a basis for the degree $n$ polynomials.
\item{2.} Interpolate the scaled polynomial $S_{1/n} P_n: \; x \mapsto P_n (\frac1n\ x)$,  at the integers $0,1,\cdots,n$, and evaluate the result at $nx$:
\begin{eqnarray*}
I_n[f](x) & = & I [S_{1/n} P_n] (nx) \\
& = & \sum_{i=0}^n c_i(S_{1/n} f) I[ t_i](nx).
\end{eqnarray*}
\end{description}
With the basis of monomials this becomes
\begin{eqnarray*}
I_n[f](x) & = & \sum_{j \in \ZZ} \sum_{i=0}^n c_i(S_{1/n} f) Q_i(j) \psi(j-nx),
\end{eqnarray*}
where $Q_i$, $i=0,\cdots,n$ are as given in (\ref{q1}) and (\ref{q2}). The coefficients in the expression must be interpreted as the appropriate ones for the monomial basis.

\section{Conclusion}

As main results of this paper, we have shown that the interpolant to a polynomial using a suitable kernel has polynomial coefficients. More importantly, a  kernel interpolant to a polynomial is constructible recursively, as is the way in which we express the error between the polynomial and its kernel interpolant. For the Gaussian kernel, we provide closed formulas for the coefficients of the kernel interpolant to a polynomial. These are given in terms of a new class of polynomials that closely resemble the classical probabilistic Hermite polynomials.
Via interpolating polynomials,  we find a way to construct a kernel interpolant to a function defined on an equally-spaced grid of a compact interval. In theory, this interpolant  uses shifts of the kernel on a full infinite grid. In numerical implementation, however, only a small number of shifts of the kernel centered outside of the interpolation interval is needed thanks to the rapid decay of the kernel. These have cleared the way for our future work in which we will investigate numerical aspects of this process. Our goals are to obtain stable and efficient algorithms for the computation of the interpolant, and to develop error estimates for $C^k-$functions having polynomial growth. With the stationary interpolation scheme, the error will not go to zero as the grid spacing contracts, but the errors estimate will be useful for analysing the residual approximation algorithm that is detailed in \cite{mlski}.

\bigskip

\noindent {\bf \large Appendix: Functions of exponential type}

\bigskip

In this appendix we show that we can extend the results of Section~\ref{general} beyond the polynomial case to that of functions of exponential type. For $f \in L(\RR)$, we use the following Fourier transform pair:
\begin{equation}\label{transform-pair}
\widehat{f}(\xi)=\int_{\RR} e^{-2 \pi i \xi x} f(x) dx, \quad \Check{f}(x)=\int_{\RR} e^{2 \pi i \xi x} f(\xi) d\xi.
\end{equation}
We assume that both Fourier transform and inverse Fourier transform have been properly extended to the Schwartz class of tempered distributions.
Let $\W$ denote the collection of all functions $\psi \in \con (\RR)$ that decays rapidly. That is, there exists a constant $C>0$, such that for
any $N \in \NN$, the following inequality holds true:
\begin{equation}\label{decay-rate}
|\psi (x)| \le \frac{C}{1+|x|^{N}}, \quad x \in \RR.
\end{equation}

Each $\psi \in \W$ induces a periodic function $\tilde{\psi}$ on $\RR$:
\begin{equation} \label{periodic}
\tilde{\psi} (x):=\sum_{z \in \ZZ} \psi(z) e^{ 2 \pi iz x}, \quad x \in \RR.
\end{equation}
The period of the above function is $1$. We will use $[0,1]$ as the fundamental interval. We are particularly interested in the
 subset $\W^*$ of $\W$ defined by
\begin{equation} \label{wiener}
\W^*:= \{\psi \in \W: \tilde{\psi} (x) \ne 0, \; x \in [0,1]\}.
\end{equation}

\begin{lemma}\label{key-lemma}
For each $\psi \in \W^*$, there exists a sequence of complex numbers $a_z,\; z \in \ZZ$, such that $\sum_{z \in \ZZ} |a_z z^k| < \infty$ for any $k \in \ZZ_+$, and
\[
\frac{1}{\tilde{\psi} (x)}=\sum_{z \in \ZZ} a_z e^{ 2 \pi iz x}, \quad x \in \RR.
\]
\end{lemma}

\proof
First off, Wiener's lemma \cite[p. 228]{katznelson} asserts that there exists a sequence of complex numbers $a_z,\; z \in \ZZ$, such that $\sum_{z \in \ZZ} |a_z| < \infty$, and
\[
\frac{1}{\tilde{\psi} (x)}=\sum_{z \in \ZZ} a_z e^{ 2 \pi iz x}, \quad x \in \RR.
\]
The rapid decay of the function $\psi$ ensures that the periodic function $\tilde{\psi}$ (see \eqref{periodic}) is infinitely differentiable on $\RR$. Since
$\tilde{\psi}(x) \ne 0$ for all $x\in \RR$, this property of smoothness of the function $\tilde{\psi}$ passes on to the function $\left(\tilde{\psi}\right)^{-1}$, whose Fourier coefficients $a_z, \; z \in \ZZ$,  therefore enjoy the desired decay condition.
\eop

We will refer to the inequality as displayed in \eqref{wiener} Wiener's condition.
Let $\mathcal{S}$ and $\mathcal{S}'$ denote, respectively, the Schwartz classes of functions and tempered distributions. For each given  $0 \le \sigma < \infty,$ let $\E_\sigma$ denote the class of analytic functions of exponential type $\sigma$. We will focus on a subclass $\E^*_\sigma$ of $\E_\sigma$ defined by:
\[
\E^*_\sigma:=\{f \in \E_\sigma: \exists C>0 \;\mbox{and} \; N \in \NN \; \mbox{such that} \; |f(x)| \le C (1+|x|^N), \; x \in \RR\}.
\]
\begin{pro} \label{existence}
Let $\psi \in \W^*$. For each $\sigma \ge 0$, and every $f \in \E^*_\sigma$ there exists a $g \in \E^*_\sigma$, such that
\[
f(j)=\sum_{z \in \ZZ} g(z) \psi (j-z), \quad j \in \ZZ.
\]
\end{pro}

\proof
By the Paley-Wiener theorem \cite[p. 162]{yosida}, we may write
\[
f=\widehat{T}, \; T \in \mathcal{S}', \mbox{supp} (T) \subset [-\sigma/(2 \pi), \sigma/(2 \pi)].\footnote{The factor $(2 \pi)^{-1}$ is the result of the particular format of the Fourier transform pair we use in the present paper; see Equation \eqref{transform-pair}.}
\]
Since $\psi \in \W^*$, we have
\[
\frac{1}{\tilde{\psi}(x)}=\sum_{z \in \ZZ} a_z e^{ 2 \pi iz x}, \quad x \in \RR,
 \]
in which the Fourier coefficients $a_z, \; z \in \ZZ,$ decay rapidly thanks to Lemma \ref{key-lemma}. That is, for any $k \in \ZZ_+$, we have $\sum_{z \in \ZZ} |a_z z^k| < \infty$. Thus, the following equation defines $\frac{T}{\tilde{\psi}} $ as  a Schwartz class distribution:
\[
\langle \frac{T}{\tilde{\psi}}, \phi \rangle:= \sum_{z \in \ZZ} a_z \langle T, e_z \cdot \phi \rangle, \quad \phi \in
\mathcal{S}.
\]
Here $e_z$ denotes the function $x \mapsto e^{2 \pi i zx}$. We also have
\[
\mbox{supp} \left(\frac{T}{\tilde{\psi}} \right) \subset [-\sigma/(2\pi), \sigma/(2\pi)].
\]
We use the same Paley-Wiener Theorem mentioned above to conclude that $g:=\left(\frac{T}{\tilde{\psi}} \right)^{\widehat{}}$ is an
 element of $\E^*_\sigma$. Thus, there exists a constant $C>0$ and an $N \in \NN$ such that
 \[
 g(x) \le C(1 + |x|^N), \quad x \in \RR.
 \]
  Now consider the function
\[
f^*(x):=\sum_{z \in \ZZ} g(x-z) \psi (z).
\]
Fix each fixed $M > 0$ and every $x \in [-M, M]$, we have
\begin{align*}
      & \left| \sum_{z \in \ZZ} g(x-z) \psi (z) \right| \\
  \le & C \sum_{z \in \ZZ}   \left( 1 +  |x-z|^N \right) \left( 1 +  |z|^{(N+2)} \right)^{-1}\\
\le & C_N M^N \sum_{z \in \ZZ}  \left( 1 +  |z| \right)^{-2}.
\end{align*}
Here $C_N$ is a constant depending only on $N$. Thus the series converges uniformly on every compact subset of $\RR$. Therefore the function $f^*$ is continuous on $\RR$ and has at most polynomial growth.
We calculate  its (distributional) inverse Fourier transform:
\[
(f^*)^{\vee} = (g)^{\vee} \cdot \tilde{\psi} = \frac{T}{\tilde{\psi}}\cdot \tilde{\psi}
= T.
\]
This shows that both $f$ and $f^*$ are the Fourier transform of the distribution $T$, meaning that they are the same function. In particular, we have
\[
\sum_{z \in \ZZ} g(j-z) \psi (z) = \sum_{z \in \ZZ} g(z) \psi (j-z) = f(j), \quad j \in \ZZ.
\]
This completes the proof.
\eop

For the uniqueness of the coefficients, we have the following result.
\begin{pro} \label{uniqueness}
Assume that $0<\epsilon < \pi$.
Let  $g \in \E^*_{\pi- \epsilon}$, $\psi \in \W^*$, and let $f$ be defined by
\[
f(x):=\sum_{z \in \ZZ} g(z) \psi (x-z),
\]
Then, in order that $f(j)=0, \; j \in \ZZ$, it is necessary and sufficient that $g(\zeta)=0, \; \zeta \in \CC$.
\end{pro}

\proof
Of course, only the necessity part needs a proof. Assume that $f(j)=0, \; j \in \ZZ$. Write
\[
f(j)=\sum_{z \in \ZZ} g(z) \psi (j-z)= \sum_{z \in \ZZ} g(j-z) \psi (z), \; j \in \ZZ,
\]
and consider the function $f^*$ defined by,
\[
f^*(x)=\sum_{z \in \ZZ} g(x-z) \psi (z).
\]
We remind readers that $f(j)=f^*(j),\; j \in \ZZ$. From the proof of Proposition \ref{existence}, we observe that $f^*$ is continuous on $\RR$ and has at most polynomial growth.  The (distributional) Fourier transform of $f^*$ can be easily calculated to be
\[
\widehat{f^*}(\xi):=\widehat{g }(\xi)\sum_{z \in \ZZ}  \psi (z) e^{ 2 \pi iz x} =\widehat{g }(\xi) \tilde{\psi}(\xi).
\]
Since $g \in E^*_{\pi-\epsilon}$, $\hat g$ is supported in $[-1/2 +\epsilon, 1/2-\epsilon]$. Thus $f^* \in E^*_{\pi-\epsilon}$. Resorting to Carlson's theorem,\footnote{Carlson's theorem asserts that a function in $E_{\pi-\epsilon}$ that vanishes on all the positive integers is identically zero; see \cite{rubel}.} we conclude that $f^*(\zeta)=0, \; \zeta \in \CC.$ The Fourier transform of $f^*$ is therefore also zero. Since
 $\widehat{f^*}(\zeta) =\widehat{g }(\zeta) \tilde{\psi}(\zeta)$, and $\tilde{\psi}(\zeta) \ne 0$, we have $\widehat{g }(\zeta)=0, \;  \zeta \in \RR$. That is, $\widehat{g}$ is the zero distribution. Hence $g$ is identically zero.
\eop

Propositions \ref{existence} and \ref{uniqueness} imply the following result.
\begin{corollary} \label{existence-and-uniqueness}
Let $0<\epsilon < \pi$ be given, and let $\psi \in \W^*$. For each $g \in \E^*_{\pi-\epsilon}$, there exists a unique $f \in \E^*_{\pi-\epsilon}$, such that
\[
\sum_{z \in \ZZ} \psi(j - z) f(z)=  g(j), \quad j \in \ZZ.
\]
\end{corollary}

Suppose that $f$ is radial (even), and that for some $\delta >0$ we have
\[
|f(x)| \le A (1+x^2)^{-(1/2 + \delta)}, \quad |\widehat{f}(\xi)| \le A (1+\xi^2)^{-(1/2 + \delta)},
\]
where $A>0$ is a constant. Then the following Poisson summation formula holds true; see \cite[p.252]{stein}.
\[
\sum_{z \in \ZZ} f(z)e^{2 \pi i z x}=\sum_{z \in \ZZ}\widehat{f}(x+z).
\]
Thus, Wiener's condition is satisfied if both $\psi$ and $\widehat{\psi}$ have the decay rate shown in \eqref{decay-rate}, and $\widehat{\psi}$ is positive. Specifically, the Gaussian kernel satisfies this condition.

We remind readers that for any $0<\epsilon < \pi$, and any polynomial $p$, we have $p \in \E^*_{\pi-\epsilon}$. To interpolate a polynomial, we do not need any extra decay condition other than what has been imposed on functions from $\W^*$. Moreover, any $\psi \in {\cal S}$ satisfying Wiener's condition suffices.

\end{document}